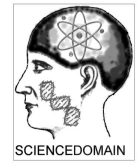

# Weak Moment of a Class of Stochastic Heat Equation with Martingale-valued Harmonic Function

**Ejighikeme Mcsylvester Omaba[1*]**

[1]*Department of Mathematics, Computer Science, Statistics and Informatics, Faculty of Science, Federal University, Ndufu-Alike Ikwo, P.M.B 1010, Abakaliki, Ebonyi State, Nigeria.*

***Author's contribution***

*The sole author designed, analyzed and interpreted and prepared the manuscript.*

***Article Information***

___

## Abstract

A study of a non-linear parabolic SPDEs of the form $\partial_t u = \mathsf{L} u + \sigma(u) f(B_t^x, t)\dot{w}$ with $\dot{w}$ as the space-time white noise and $f(B_t^x, t)$ a space-time harmonic function was done. The function $\sigma : \mathbf{R} \to \mathbf{R}$ is Lipschitz continuous and $\mathsf{L}$ the $L^2$-generator of a Lévy process. Some precise condition for existence and uniqueness of the solution were given and we show that the solution grows weakly(in law/distribution) in time (for large $t$) at most a precise exponential rate for the $\mathsf{L}$; and grows in time at most a precise exponential rate for the case of $\mathsf{L} = -(-\Delta)^{\alpha/2}, \alpha \in (1,2]$ generator of an alpha-stable process.

# 1 Introduction

The study of stochastic heat equations with space-time white noise has received quite a lot of interests, see [1-4] and [3-8], to mention but a few. Here, we study a long-time behaviour of a class of stochastic heat equation perturbed with a multiplicative space-time harmonic function of $f(B_t^x, t)$ type. The space-time harmonic functions are closely related to a class of martingales that are instantaneous functions of $X_t$. That is, every space-time harmonic function $f$ gives rise to a martingale with respect to a probability measure. Thus given a process $X_t$, $f(X_t, t)$ is a martingale with respect to a probability measure $d\mu(\lambda)$. Space-time harmonic functions have their application in determining the growth properties of a branching process [9]. They also play a central role in constructing Schrödinger bridges for Markov chains. An example of the space-time harmonic function is the Brownian bridge given by $f(x,t) := \sqrt{1-t} \exp(x^2/2(1-t))$. Martingale on the other hand plays a crucial role in the Black-Scholes formula. The continuous-time stochastic process $S_t$ given below describes the price of the stock (stock price process) at time $t$, $S_0$ the current stock price, known to all investors at $t = 0$:

$$S_t = S_0 \exp[(\mu - \lambda^2/2)t + \lambda B_t], \ t \in [0,T]$$

or equivalently

$$\log \frac{S_t}{S_0} = (\mu - \lambda^2/2)t + \lambda B_t, \ t \in [0,T]$$

the "log-normal" or "geometric Brownian motion" where $B_t$ denotes a Brownian motion, the mean $\mu$ denotes the "drift term" and $\lambda$ the volatility. The stock price $S_t$ can be thought of as an exponential of Brownian motion with drift term. That is, the geometric Brownian motion models a stock price. It has the following martingale property (known as its expected growth):

$$ES_t = S_0 e^{\mu t}, \ for \ t \geq 0, \ since \ Ee^{\lambda B_t} = e^{\lambda^2 t/2}.$$

In [10], the authors studied space-time harmonic functions of both Brownian and Gamma-types and their applications to finance. Due to the importance and applications of the martingale-valued harmonic function(processes), for example, Variance-Gamma (VG) process is used as a model for asset returns (log-price increments) and options pricing; we are motivated to investigate a long-time influence and effect of the space-time harmonic functions for Brownian motion on a class of stochastic heat equation.

Now consider a parabolic heat equation perturbed with a multiplicative space-time harmonic function.

$$\partial_t u(x,t) = \mathsf{L} u(x,t) + \sigma(u(x,t)) f(B_t^x, t) \dot{w}(t,x), \tag{1.1}$$

with the non-random initial condition $u(x,0) = u_0(x)$, $\sigma : \mathbf{R} \to \mathbf{R}$ a Lipschitz continuous function and

$$f(B_t^x, t) = \int_{\mathbf{R}} M_t^\lambda d\mu(\lambda).$$

Here, $B_t^x$ does not stand for a Brownian motion starting from $x$, rather a notation specifying space dependence $x$. Most times, we write $f(B_t, t)$ without the space variable $x$. The aim of this article is to study how the solution $u$ is influenced by the interaction or competition between the generator of a semigroup $L$ which has a smoothing effect, the martingale measure and the noise potential $w(dxdt)$, which makes the solution spatially irregular. This is a long-time behaviour of a system exhibiting an intermittency effect and we show that existence of the solution is closely related to its growth.

**Definition 1.1:** *[10] Let $X = (X_t, t \geq 0)$ be a real-valued stochastic process, we say that $f : \mathbf{R} \times \mathbf{R}_+ \to \mathbf{R}$ is a space-time harmonic function if $\{f(X_t, t), t \geq 0\}$ is a martingale (with respect to the filtration of $X$).*

Here, follows our representation of the harmonic function (see [10] for details on space-time harmonic functions):

**Theorem 1.2:** *Every $\mathbf{R}_+$-valued space-time harmonic function of $(B_t, t \geq 0)$ such that $f(0,0) = 1$, may be written as*

$$f(B_t, t) = \int_{\mathbf{R}} d\mu(\lambda) \exp(\lambda B_t - \frac{\lambda^2}{2} t),$$

*for some probability measure $d\mu(\lambda)$ on $\mathbf{R}$.*

The process $M_t^\lambda = \exp(\lambda B_t - \frac{\lambda^2}{2} t)$ known as an exponential martingale has its natural appearance in the expression of the Randon-Nikodym densities with the following property:

$$E[M_t^\lambda] = E[M_0^\lambda] = 1, \forall\, t \geq 0.$$

One approach to solving the problem is to express $M_t^\lambda$ in the following classical series expansion

$$\exp(\lambda B_t - \frac{\lambda^2}{2} t) = \sum_{n=0}^{\infty} \frac{\lambda^n}{n!} H_n(B_t, t)$$

where $H_n(B_t, t) = t^{n/2} h_n(\frac{B_t}{\sqrt{t}})$ with $(h_n)$ the Hermite polynomials and $H_n(B_t, t)$ a martingale and seek for an approximation of the hermite polynomial.

**Definition 1.3:** *Hermite polynomials $(h_n)$ on $\mathbf{R}$ is given by*

$$h_n(x) = (-1)^n e^{x^2/2} \frac{d^n}{dx^n} e^{-x^2/2}, \quad \text{for } n = 0, 1, \ldots.$$

The first few terms are

$$h_0(x) = 1, \ h_1(x) = x, \ h_2(x) = x^2 - 1, \ h_3(x) = x^3 - 3x, \ etc$$

with the following recursive formula $h_{n+1}(x) = xh_n(x) - nh_{n-1}(x)$ and the generating function $e^{xt-t^2/2} = \sum_{n=0}^{\infty} \frac{t^n}{n!} h_n(x)$. Alternatively, one gives the generating identity of the Hermite functions $h_n(x)$ as follows $e^{-(x^2/2-tx-t^2)} = \sum_{n=0}^{\infty} \frac{t^n}{n!} h_n(x)$. They satisfy the so called "creation" and "annihilation" identities

$$(x - \frac{d}{dx})h_n(x) = h_{n+1}(x) \ \ and \ \ (x + \frac{d}{dx})h_n(x) = h_{n-1}(x), \ \forall \ n \geq 0$$

where $h_{-1}(x) = 0$, $h_0(x) = e^{-x^2/2}$, $h_1(x) = 2xe^{-x^2/2}$, etc. Generally,

$$h_n(x) = P_n(x)e^{-x^2/2},$$

where $P_n$ is a polynomial of degree $n$. The existence and uniqueness of the solution requires a bound on the polynomial function $P_n$ for all $n$. We will follow rather a different approach in this paper.

The outline of the paper is given as follows. Section 2 is the formulation of the problem and the main results of the paper, section 3 surveys some preliminary concepts used. Auxiliary results which comprise of Lemma(s) and Propositions that will be used in the proofs of the results are given in section 4. Proofs of main results: existence and uniqueness result and growth bounds are given in the last section.

## 2 Formulation and Main Results

For existence and uniqueness, we need the following condition on $\sigma$. The condition states that $\sigma$ is globally Lipschitz in its variable.

**Condition 2.1:** *There exist a finite positive constant,* $\text{Lip}_\sigma$ *such that for all* $x, y \in \mathbf{R}$, *we have*

$$|\sigma(x) - \sigma(y)| \leq \text{Lip}_\sigma |x - y|. \tag{2.1}$$

We will make use of Walsh's integral in [2] and follow the steps of [1].

**Definition 2.2:** *We say that a process* $\{u(t,x)\}_{x \in \mathbf{R}, t>0}$ *is a mild solution of (1.1) if a.s, the following is satisfied*

$$u(x,t) = (P_t u_0)(x) + \int_0^t \int_{\mathbf{R}} p(t-s, x, y) \sigma(u(y,s)) f(B_s^y, s) w(dy ds) \tag{2.2}$$

where $p(t,.,.)$ is the heat kernel. If in addition to the above, $\{u(x,t)\}_{x \in \mathbf{R}, t>0}$ satisfies the following condition

$$\sup_{0 \leq t \leq T} \sup_{x \in \mathbf{R}} \mathrm{E} \, |u(x,t)|^2 < \infty, \tag{2.3}$$

for all $T > 0$, then we say that $\{u(x,t)\}_{x \in \mathbf{R}, t > 0}$ is a random field solution to (1.1).

Define

$$\Upsilon(\beta) := \frac{1}{2\pi} \int_{\mathbf{R}} \frac{\mathrm{d}\xi}{\beta + 2\mathbf{Re}\Psi(\xi)} \quad \text{for all } \beta > 0, \tag{2.4}$$

where $\Psi$ is the characteristic exponent for the Lévy process. A result of Dalang [3] shows that equation (1.1) has a unique solution with the requirement that $\Upsilon(\beta) < \infty$ for all $\beta > 0$. Fix some $x_0 \in \mathbf{R}$ and define the *upper $p$ th-moment Liapunov exponent* $\bar{\gamma}(p)$ of $u$ [at $x_0$] as

$$\bar{\gamma}(p) := \limsup_{t \to \infty} \frac{1}{t} \ln \mathrm{E}\left[|u(x_0,t)|^p\right] \quad \text{for all } p \in (0, \infty). \tag{2.5}$$

See [1] and its references for details on intermittent property. We will assume the following condition:

$$f(B_s, s) = \int_{\mathbf{R}} M_s^\lambda \mathrm{d}\mu(\lambda) \leq \sup_{\lambda \in \mathbf{R}} M_s^\lambda := M_s^{\lambda_0}$$

for a fixed $\lambda_0$ finite. Consequently,

$$u(x,t) \leq (P_t u_0)(x) + \sup_{0 < s < t} M_s^{\lambda_0} \int_0^t \int_{\mathbf{R}} p(t-s, x, y) \sigma(u(y,s)) w(\mathrm{d}y \mathrm{d}s).$$

**Theorem 2.3:** *Equation (1.1) admits a mild solution $u$ that is unique up to a modification, satisfying the following: For even integers $p > 1$,*

$$\bar{\gamma}(p) \leq \inf\{\beta > 0 : \Upsilon(\frac{2\beta}{p}) < \frac{1}{(z_p \, C \, \mathrm{Lip}_\sigma(\frac{p}{p-1}) \exp(\lambda_0^2 t_0 \, (p-1)/2))^2}\},$$

where $z_p$ is the constant given in [1]. Next, we estimate bound on growth moment and show that our solution grows exponentially. For the lower bound result, we will need the following extra condition on $\sigma$.

**Condition 2.4:** *There exist a positive constant, $L_\sigma$ such that for all $x \in \mathbf{R}$, we have*

$$|\sigma(x)| \geq L_\sigma |x| \tag{2.6}$$

**Theorem 2.5:** *If we further assume that condition 2.4 holds and $\inf_{x \in \mathbf{R}} u_0(x) > 0$, then in law (in distribution) we have:*

$$\overline{\gamma}(2) \geq \Upsilon^{-1}\left(\frac{1}{\mathsf{K}^2}\right) > 0,$$

where $\Upsilon^{-1}(t) := \sup\{\beta > 0 : \Upsilon(\beta) > t\}$ and $\mathsf{K}^2 = L_\sigma^2 C^2 \sqrt{\pi} \exp\left(a - \frac{|a|}{\sqrt{2}}\right)^2$.

To prove the lower bound for the solution, we consider a heat equation of the form

$$\partial_t u(x,t) = \mathsf{L} u(x,t) + \sigma(u(x,t)) \dot{w}(t,x) \int_{\mathbf{R}} d\mu(\lambda) \exp\left(\lambda B_t - \frac{\lambda^2}{2} t\right), \qquad (2.7)$$

with the non-random initial condition $u(x,0) = u_0(x)$. For explicit and easy computation, we will take $d\mu(\lambda)$, the probability measure of a normalized (standardized) Gaussian random variable:

$$d\mu(\lambda) = C \exp(-\lambda^2/2) d\lambda, \ C > 0.$$

The solution to equation (2.7) is given by

$$\begin{aligned} u(x,t) &= \int_{\mathbf{R}} p(t,x,y) u_0(y) dy \\ &\quad + C \int_0^t \int_{\mathbf{R}} \int_{\mathbf{R}} p(t-s,x,y) \sigma(u(y,s)) \exp\left(\lambda B_s - \frac{\lambda^2}{2}(1+s)\right) d\lambda w(dyds). \end{aligned}$$

The approach we will adopt to making sense of our solution is in terms of a stopping time, particularly (and more importantly in terms of applications) a hitting time of a Lévy process. Let a new process $X = X_T(t), \ t \geq 0$ be defined by $X_T(t) = X_{T+t} - X_T$, (see [11]).

**Theorem 2.6:** *(Strong Markov property) If $X$ is a Lévy process and $T$ is a stopping time, then on $(T < \infty)$,*

1. $X_T$ is a Lévy process that is independent of $\mathsf{F}_T$,
2. For each $t \geq 0$, $X_T(t) = X_{T \wedge t}$ has the same law as $X_t$,
3. $X_T$ has cádlág paths and is $\mathsf{F}_{T+t}$-adapted.

**Remark 2.7:** *Stating the above theorem for the hitting time $T_a$ of level "$a$" of a process $B_t$.*

1. For each $t \geq 0$, $B_{T \wedge t}$ has the same law as $B_t$.

**Definition 2.8:** *(Hitting Times) Let $a \in \mathbf{R}$ and define*

$$T_a := \inf\{t \geq 0 : B_t = a\}.$$

The random variable $T_a$ is called the first hitting time of level "$a$" by a Brownian motion.

**Proposition 2.9:** *[12,13] Let $T_a$ be the first hitting time. Then $P(T_a < \infty) = 1$ and $\mathrm{E}[T_a] = \infty$. Moreover,*

$$\mathrm{E}[\exp(-\lambda T_a)] = \exp(-|a|\sqrt{2\lambda}),\ \lambda \geq 0.$$

**Remark 2.10:** *On the set $\{T_a < \infty\}$, we have that $B_{t \wedge T_a} \to B_{T_a} = a$ when $t \to \infty$, hence*

$$\exp(\lambda B_t - \frac{\lambda^2}{2}t) \to \exp(\lambda a - \frac{\lambda^2}{2}T_a).$$

In what follows, we consider the case of $\mathsf{L} = -(-\Delta)^{\alpha/2}$, $\alpha \in (1,2]$ a generator of an alpha-stable Lévy processes and use some explicit bounds on its corresponding fractional heat kernel to obtain more precise results. It helps to understand the full behaviour of the solution and to have an explicit estimate for the generator of the process. We present some required properties of $p(t,x)$ which come in handy in the proof of our results [14, 15].

**Lemma 2.11:** *[15] Suppose that $p(t,x)$ denotes the heat kernel for a strictly stable process of order $\alpha$. Then the following estimate holds.*

$$p(t,x,y) \approx t^{-d/\alpha} \wedge \frac{t}{|x-y|^{d+\alpha}} \quad \text{for all} \quad t > 0 \quad \text{and} \quad x, y \in \mathbf{R}^d.$$

Here and in the sequel, for two non-negative functions $f$, $g$, $f \approx g$ means that there exists a positive constant $c > 1$ such that $c^{-1}g \leq f \leq cg$ on their common domain of definition. Below is the upper growth bound of the solution.

**Theorem 2.12:** *There exist constants $c$ and $c'$ such that*

$$\sup_{x \in \mathbf{R}} \mathrm{E}|u(x,t)|^p \leq c \exp(c'(\frac{p}{p-1})^{2\alpha/(\alpha-1)}t) \text{ for all } t > 0,$$

where the constant $c'$ depends on $p$, $\lambda_0$ and $t_0$.

## 3 Preliminaries

Let $(\Omega, \mathsf{F}, \mathsf{P})$ be a probability space. Let $B_t$ be a real-valued random process with index set $\mathbf{R}_+ = [0, \infty)$ and consider $\Omega := C(\mathbf{R}_+, \mathbf{R})$, the space of real-valued continuous functions on $\mathbf{R}_+$, equipped with the following metric

$$\rho(f,g) := \sum_{k \geq 0} \frac{1}{2^k} \sup_{0 \leq |x| \leq k} |f(x) - g(x)| \wedge 1.$$

The filtration $\mathsf{F}_t$ is given as a subset of the algebra, thus: $\mathsf{F} := \mathsf{B}(\Omega)$. Brownian motion is exponentially integrable with

$$\psi(\xi) := \mathrm{E}[e^{i\xi B_t}] = e^{-\xi^2 t/2}.$$

Martingale is a very important concept in stochastic analysis because they describe fair games. It is a process where the current state $X_t$ gives the best prediction for its future state. We state one important fact about martingales. See [11, 12, 13] for details.

**Lemma 3.1:** *(Doob's Maximal Inequality) Suppose $\{X_t, t \geq 0\}$ is a continuous sub-martingale and $p > 1$. Then for any $t \geq 0$,*

$$\mathrm{E}[\sup_{0 < s < t} X_s]^p \leq (\frac{p}{p-1})^p \mathrm{E}[X_t^p], \ p > 1$$

provided $X_s \geq 0$, a.s $P$ for every $s \geq 0$ and $\mathrm{E}[X_t^p] < \infty$.

Let $S_t = \sup_{s<t} B_t$.

**Proposition 3.2:** *Then for $a > 0$,*

$$\mathrm{P}[S_t \geq at] \leq \exp(-a^2 t/2).$$

**Remark 3.3:**

We note the following identities:

$$\exp(\lambda S_t - \lambda^2 t/2) \leq \sup_{s \leq t} M_s^\lambda \ \ and \ \ \inf_{\lambda > 0}(-\lambda at + \lambda^2 t/2) = -a^2 t/2.$$

We now present the following renewal inequality from ([16], chapter 7). Here, we desire bound on the function involved rather than finding their asymptotic properties.

**Proposition 3.4:** *[4] Let $\rho > 0$ and suppose $f(t)$ is a non-negative and locally integrable function satisfying*

$$f(t) \leq c_1 + \kappa \int_0^t (t-s)^{\rho-1} f(s) \mathrm{d}s \ \ for \ all \ t > 0,$$

*where $c_1$ is some positive number. Then we have*

$$f(t) \leq c_2 \exp(c_3 (\Gamma(\rho))^{1/\rho} \kappa^{1/\rho} t) \ for \ all \ t > 0,$$

*for some positive constants $c_2$ and $c_3$.*

## 4 Auxiliary Results

Define for $t > 0$,

$$(\mathsf{A}u)(x,t) = \sup_{0<s<t} M_s^{\lambda_0} \int_0^t \int_{\mathbf{R}} p(t-s,x,y)\sigma(u(y,s))w(\mathrm{d}y\mathrm{d}s).$$

We will also use the following norm:

$$\|u\|_{p,\beta} := \left\{\sup_{t>0}\sup_{x\in\mathbf{R}} e^{-\beta t} \mathrm{E}[|u(x,t)|^p]\right\}^{1/p} \text{ for } p \geq 2. \tag{4.1}$$

**Lemma 4.1:** *For all* $\beta > 0$,

$$\sup_{t>0} e^{-\beta t} \int_0^t \int_{\mathbf{R}} |p(s,x,y)|^2 \, \mathrm{d}y\mathrm{d}s \leq \Upsilon(\beta).$$

*Proof.*

$$\sup_{t>0} e^{-\beta t} \int_0^t \int_{\mathbf{R}} |p(s,x,y)|^2 \, \mathrm{d}y\mathrm{d}s \leq \int_0^\infty \int_{\mathbf{R}} e^{-\beta s} |p(s,x,y)|^2 \, \mathrm{d}y\mathrm{d}s$$

$$= \int_0^\infty \int_{\mathbf{R}} e^{-\beta s} |\hat{p}(s,\xi)|^2 \, \mathrm{d}\xi\mathrm{d}s$$

$$= \int_{\mathbf{R}} \frac{\mathrm{d}\xi}{\beta + 2\mathrm{Re}\Psi(\xi)} = \Upsilon(\beta).$$

**Proposition 4.2:** *Suppose that* $u$ *admits a predictable version and that* $\|u\|_{p,\beta} < \infty$ *for* $\beta > 0$. *Then there exists some positive constant* $z_p$ *such that*

$$\|\mathsf{A}u\|_{p,\beta} \leq z_p(\frac{p}{p-1})\exp(\lambda_0^2 t_0(p-1)/2)[\mathrm{Lip}_\sigma \|u\|_{p,\beta} + |\sigma(0)|]\sqrt{\Upsilon(2\beta/p)}.$$

*Proof.*

$$\mathrm{E}|(\mathsf{A}u)(x,t)|^p \leq z_p^p(\frac{p}{p-1})^p \exp(\frac{p\lambda_0^2 t}{2}(p-1))$$

$$\times (\int_0^t \int_{\mathbf{R}} |p(t-s,x,y)|^2 [\mathrm{E}|\sigma(u(s,y))|^p]^{2/p} \, \mathrm{d}y\mathrm{d}s)^{p/2}.$$

Therefore,

$$[\mathrm{E}\,|\,(\mathrm{A}u)(x,t)\,|^p\,]^{2/p}$$
$$\leq z_p^2 \left(\frac{p}{p-1}\right)^2 \exp(\lambda_0^2 t(p-1)) \int_0^t \int_{\mathbf{R}} |\,p(t-s,x,y)\,|^2 \,[\mathrm{E}\,|\,\sigma(u(y,s))\,|^p\,]^{2/p}\,\mathrm{d}y\mathrm{d}s.$$

$$\leq z_p^2 \left(\frac{p}{p-1}\right)^2 \exp(\lambda_0^2 t(p-1))$$
$$\times \int_0^t \int_{\mathbf{R}} |\,p(t-s,x,y)\,|^2 \,[\mathrm{Lip}_\sigma [\mathrm{E}\,|\,u(y,s)\,|^p\,]^{1/p} + |\,\sigma(0)\,|]^2 \,\mathrm{d}y\mathrm{d}s$$
$$\leq z_p^2 \left(\frac{p}{p-1}\right)^2 \exp(\lambda_0^2 t(p-1))$$
$$\times \int_0^t \int_{\mathbf{R}} |\,p(t-s,x,y)\,|^2 \,[\mathrm{Lip}_\sigma^2 [\mathrm{E}\,|\,u(y,s)\,|^p\,]^{2/p} + |\,\sigma(0)\,|^2]\,\mathrm{d}y\mathrm{d}s.$$

Multiply through by $\exp(-2\beta t/p)$,

$$[\mathrm{e}^{-\beta t}\mathrm{E}\,|\,(\mathrm{A}u)(x,t)\,|^p\,]^{2/p}$$
$$\leq z_p^2 \left(\frac{p}{p-1}\right)^2 \exp(\lambda_0^2 t(p-1))\mathrm{e}^{-2\beta/p(t-s)}$$
$$\times \int_0^t \int_{\mathbf{R}} |\,p(t-s,x,y)\,|^2 \,\mathrm{e}^{-2\beta/ps}[\mathrm{Lip}_\sigma^2 [\mathrm{E}\,|\,u(y,s)\,|^p\,]^{2/p} + |\,\sigma(0)\,|^2]\,\mathrm{d}y\mathrm{d}s.$$

It follows that for fixed $t_0 \geq 0$

$$\|\,\mathrm{A}u\,\|_{p,\beta}^2 \leq z_p^2 \left(\frac{p}{p-1}\right)^2 \exp(\lambda_0^2 t_0(p-1))[\mathrm{Lip}_\sigma^2 \,\|\,u\,\|_{p,\beta}^2 + |\,\sigma(0)\,|^2]$$
$$\times \sup_{x\in\mathbf{R}} \sup_{t\geq 0} \int_0^t \int_{\mathbf{R}} \mathrm{e}^{-2\beta/p(t-s)} |\,p(t-s,x,y)\,|^2 \,\mathrm{d}y\mathrm{d}s$$
$$\leq z_p^2 \left(\frac{p}{p-1}\right)^2 \exp(\lambda_0^2 t_0(p-1))[\mathrm{Lip}_\sigma^2 \,\|\,u\,\|_{p,\beta}^2 + |\,\sigma(0)\,|^2]$$
$$\times \int_0^\infty \int_{\mathbf{R}} \mathrm{e}^{-2\beta/ps} |\,\hat{p}(s,\xi)\,|^2 \,\mathrm{d}\xi\mathrm{d}s$$
$$= z_p^2 \left(\frac{p}{p-1}\right)^2 \exp(\lambda_0^2 t_0(p-1))[\mathrm{Lip}_\sigma^2 \,\|\,u\,\|_{p,\beta}^2 + |\,\sigma(0)\,|^2]\Upsilon(2\beta/p).$$

**Proposition 4.3:** *Let $\beta > 0$ and let $u$ and $v$ be two predictable random field solutions satisfying* $\|\,u\,\|_{p,\beta} + \|\,v\,\|_{p,\beta} < \infty$. *Then*

$$\|\,\mathrm{A}u - \mathrm{A}v\,\|_{p,\beta} \leq z_p \mathrm{Lip}_\sigma \left(\frac{p}{p-1}\right) \exp(\lambda_0^2 t_0(p-1)/2)\,\|\,u-v\,\|_{p,\beta}\,\sqrt{\Upsilon(2\beta/p)}.$$

*Proof.* We continue in similar manner as above.

$$[\mathrm{E} \,|\, (\mathcal{A}u)(x,t) - \mathcal{A}v(x,t) \,|^p]^{2/p}$$
$$\leq (z_p \mathrm{Lip}_\sigma)^2 (\frac{p}{p-1})^2 \exp(p\lambda_0^2 t(p-1))$$
$$\times \int_0^t \int_{\mathbf{R}} |\, p(t-s,x,y)\,|^2 \,[\mathrm{E}\,|\, u(y,s) - v(y,s)\,|^p]^{2/p}\, \mathrm{d}y\mathrm{d}s.$$

Multiplying through by $\exp(-2\beta t/p)$, therefore

$$[\mathrm{e}^{-\beta t}\mathrm{E} \,|\, (\mathcal{A}u)(x,t) - \mathcal{A}v(x,t)\,|^p]^{2/p}$$
$$\leq (z_p \mathrm{Lip}_\sigma)^2 (\frac{p}{p-1})^2 \exp(p\lambda_0^2 t(p-1))$$
$$\times \mathrm{e}^{-2\beta(t-s)/p} \int_0^t \int_{\mathbf{R}} |\, p(t-s,x,y)\,|^2 \,[\mathrm{e}^{-\beta s}\mathrm{E}\,|\, u(y,s) - v(y,s)\,|^p]^{2/p}\, \mathrm{d}y\mathrm{d}s.$$

It follows that

$$\|\mathcal{A}u - \mathcal{A}v\|^2_{p,\beta} \leq (z_p \mathrm{Lip}_\sigma)^2 (\frac{p}{p-1})^2 \exp(p\lambda_0^2 t_0 (p-1))$$
$$\times \|u - v\|^2_{p,\beta} \int_0^t \int_{\mathbf{R}} \mathrm{e}^{-2\beta s/p} |\,\hat{p}(s,\xi)\,|^2 \, \mathrm{d}\xi\mathrm{d}s.$$

## 5 Proofs of Main Results

The proofs of the main results of this paper are outlined.

### 5.1 Existence and uniqueness

Firstly, we establish the existence and uniqueness of equation (1.1) under linear condition on $\sigma$.

*Proof of the existence and uniqueness part of Theorem 2.3.* We prove the existence of the solution by an iterative schemes. Let's define $v_0(x,t) := u_0(x)$ for all $t \geq 0$ and $x \in \mathbf{R}$. Since $u_0$ is assumed to be bounded, so $\|v_0\|_{p,\beta} < \infty$ for all $\beta > 0$. Iteratively, we set

$$\begin{cases} v_{n+1}(x,t) &= \mathcal{A}v_n(x,t) + (P_t u_0)(x) \\ v_n(x,t) &= \mathcal{A}v_{n-1}(x,t) + (P_t u_0)(x). \end{cases}$$

From the above, we have that for sufficiently large $\beta$,

$$\|\mathcal{A}v_{n+1}\|_{p,\beta} = \|\mathcal{A}(\mathcal{A}v_n)\|_{p,\beta}.$$

Then by Proposition 4.2, we have that

$$\| Av_{n+1} \|_{p,\beta} = \| A(Av_n) \|_{p,\beta} \leq z_p(\frac{p}{p-1})\exp(\lambda_0^2 t_0(p-1)/2) \quad (5.1)$$
$$\times [\mathrm{Lip}_\sigma \| Av_n \|_{p,\beta} + |\sigma(0)|]\sqrt{\Upsilon(2\beta/p)}.$$

Since $\lim_{\beta \to \infty} \Upsilon(\beta) = 0$, then we can always choose and fix $\beta > 0$ such that

$$z_p(\frac{p}{p-1})\exp(\lambda_0^2 t_0(p-1)/2)\sqrt{\Upsilon(2\beta/p)}\mathrm{Lip}_\sigma < 1 \Leftrightarrow Q_p^2 \Upsilon(2\beta/p)\mathrm{Lip}_\sigma^2 < 1$$

where $Q_p = z_p(\frac{p}{p-1})\exp(\lambda_0^2 t_0(p-1)/2)$ since $Q_p\sqrt{\Upsilon(2\beta/p)}\mathrm{Lip}_\sigma > 0$. It follows that

$$\Upsilon(2\beta/p) < \frac{1}{[Q_p \mathrm{Lip}_\sigma]^2} < \infty, \forall \beta > 0.$$

From (5.1) we have $\| Av_{n+1} \|_{p,\beta} \leq Q_p\sqrt{\Upsilon(2\beta/p)}\mathrm{Lip}_\sigma \| Av_n \|_{p,\beta} + Q_p\sqrt{\Upsilon(2\beta/p)}$. Taking sup of both sides over $n$, therefore,

$$\sup_{n \geq 0} \| Av_{n+1} \|_{p,\beta} - Q_p\sqrt{\Upsilon(2\beta/p)}\mathrm{Lip}_\sigma \sup_{n \geq 0} \| Av_n \|_{p,\beta} \leq Q_p\sqrt{\Upsilon(2\beta/p)}.$$

Rest of the proof of the existence result follows from the proof of Theorem 3.1.1 of the thesis [13]. For the proof of the uniqueness of the solution up to modification. Let $u_1$ and $u_2$ be solutions and assume for contradiction that $u_1 \neq u_2$ such that

$$\begin{cases} u_1(x,t) &= Au_1(x,t) + (P_t u_0)(x) \\ u_2(x,t) &= Au_2(x,t) + (P_t u_0)(x). \end{cases}$$

Therefore, $\| u_1 - u_2 \|_{p,\beta} = \| Au_1 - Au_2 \|_{p,\beta}$, $\beta > 0$. Then by Proposition 4.3 we have that $\| u_1 - u_2 \|_{p,\beta} \leq Q_p \mathrm{Lip}_\sigma \sqrt{\Upsilon(2\beta/p)} \| u_1 - u_2 \|_{p,\beta}$, and $\| u_1 - u_2 \|_{p,\beta}[1 - Q_p \mathrm{Lip}_\sigma \sqrt{\Upsilon(2\beta/p)}] \leq 0$. This implies that $\| u_1 - u_2 \|_{p,\beta} \leq 0$ (since $1 - Q_p \mathrm{Lip}_\sigma \sqrt{\Upsilon(2\beta/p)} > 0$) which implies that $\| u_1 - u_2 \|_{p,\beta} = 0$, and follows that $u_1 = u_2$. This contradicts the assumption that $u_1 \neq u_2$, hence $u_1 = u_2$ and thus a unique solution. Therefore $u_1$ and $u_2$ are modification of each other. The proof of the upper bound result of the theorem follows similar argument as in the proof of Theorem 2.5 below.

### 5.2 Growth bounds

Here we give bounds on the growth of the second moment of the equation.

*Proof of Theorem 2.5.* The proof suffices for us to show that

$$\int_0^\infty e^{-\beta t} \mathrm{E}\,|u(x,t)|^2\,\mathrm{d}t = \infty,\ \forall\, t>0$$

whenever $\Upsilon(\beta) \geq \dfrac{1}{\mathsf{K}^2}$ with $\mathsf{K}^2 := L_\sigma^2 C^2 \sqrt{\pi}\exp(a-\dfrac{|a|}{\sqrt{2}})^2$. We aim to establish this claim following similar steps in [1, 4]. We therefore have that for large $t$, the solution to equation (2.7) is given by

$$u(x,t) \stackrel{d}{=} (P_t u_0)(x) + C\int_0^t\!\!\int_{\mathbf{R}}\!\!\int_{\mathbf{R}} p(t-s,x,y)\sigma(u(s,y))$$
$$\times \exp(\lambda B_{s\wedge T_a} - \frac{\lambda^2}{2}(1+(s\wedge T_a)))\mathrm{d}\lambda\,w(\mathrm{d}y\mathrm{d}s)$$
$$\stackrel{d}{=} (P_t u_0)(x) + C\int_0^t\!\!\int_{\mathbf{R}}\!\!\int_{\mathbf{R}} p(t-s,x,y)\sigma(u(s,y))$$
$$\times \exp(\lambda B_{T_a} - \frac{\lambda^2}{2}(1+T_a))\mathrm{d}\lambda\,w(\mathrm{d}y\mathrm{d}s),\ as\ s\to\infty,\{T_a<\infty\}.$$

Thus in law (distribution), we have the following

$$\mathrm{E}\,|u(x,t)|^2 = |(P_t u_0)(x)|^2 + C^2\int_0^t\!\!\int_{\mathbf{R}}\!\!\int_{\mathbf{R}} p^2(t-s,x,y)\mathrm{E}\,|\sigma(u(s,y))|^2$$
$$\times \mathrm{E}\exp(2\lambda a - \lambda^2(1+T_a))\mathrm{d}\lambda\mathrm{d}y\mathrm{d}s$$
$$= |(P_t u_0)(x)|^2 + C^2\int_0^t\!\!\int_{\mathbf{R}}\!\!\int_{\mathbf{R}} p^2(t-s,x,y)\mathrm{E}\,|\sigma(u(s,y))|^2$$
$$\times \exp(\lambda(2a-|a|\sqrt{2}) - \lambda^2)\mathrm{d}\lambda\mathrm{d}y\mathrm{d}s.$$

Hence

$$\mathrm{E}\,|u(x,t)|^2 = |(P_t u_0)(x)|^2$$
$$+ C^2\sqrt{\pi}\exp(a-\frac{|a|}{\sqrt{2}})^2 \int_0^t\!\!\int_{\mathbf{R}} p^2(t-s,x,y)\mathrm{E}\,|\sigma(u(s,y))|^2\,\mathrm{d}y\mathrm{d}s.$$

Apply Laplace transforms to both sides for all $\beta>0$ and $x\in\mathbf{R}$, and using the fact that $|\sigma(u)|\geq L_\sigma|u|$,

$$\int_0^\infty e^{-\beta t}\mathrm{E}\,|u(t,x)|^2\mathrm{d}t \geq \int_0^\infty e^{-\beta t}|(P_t u_0)(x)|^2\mathrm{d}t + L_\sigma^2 C^2\sqrt{\pi}\exp(a-\frac{|a|}{\sqrt{2}})^2$$
$$\times \int_0^\infty e^{-\beta t}\int_0^t\!\!\int_{\mathbf{R}} \mathrm{E}(|u(s,y)|^2)\,|p(t-s,x-y)|^2\,\mathrm{d}y\mathrm{d}s\mathrm{d}t$$
$$= \int_0^\infty e^{-\beta t}|(P_t u_0)(x)|^2\mathrm{d}t + L_\sigma^2 C^2\sqrt{\pi}\exp(a-\frac{|a|}{\sqrt{2}})^2$$
$$\times \int_{\mathbf{R}}\mathrm{d}y\int_0^\infty e^{-\beta s}|p(s,x-y)|^2\,\mathrm{d}s\int_0^\infty e^{-\beta s}\mathrm{E}\,|u(y,s)|^2\,\mathrm{d}s.$$

*Omaba; ARJOM, 3(1): 1-16, 2017; Article no.ARJOM.31665*Then iteratively, we have

$$F_\beta(x) \geq \frac{\varepsilon^2}{\beta}\sum_{n=0}^{\infty}(\mathsf{K}^2\Upsilon(\beta))^n \quad \text{and} \quad F_\beta(x) = \infty \quad \text{whenever} \quad \mathsf{K}^2\Upsilon(\beta) \geq 1,$$

where $F_\beta(x) := \int_0^\infty e^{-\beta t} \mathrm{E}|u(x,t)|^2 \, dt$ and the result follows.

*Proof of Theorem 2.12.* Using the condition that $u_0(x) \leq c_1$, we obtain

$$\mathrm{E}|u(x,t)|^p \leq c_1^p + z_p^p \left(\frac{p}{p-1}\right)^p \exp\left(\frac{p\lambda_0^2 t}{2}(p-1)\right)$$
$$\times \left(\int_0^t \int_{\mathbf{R}} |p(t-s,x,y)|^2 [\mathrm{E}|\sigma(u(y,s))|^p]^{2/p} \, dyds\right)^{p/2}.$$

Therefore, with $\sigma(0) = 0$,

$$[\mathrm{E}|u(x,t)|^p]^{2/p} \leq c_1^p + z_p^2 \left(\frac{p}{p-1}\right)^2 \exp(\lambda_0^2 t(p-1))$$
$$\times \int_0^t \int_{\mathbf{R}} |p(t-s,x,y)|^2 [\mathrm{E}|\sigma(u(y,s))|^p]^{2/p} \, dyds$$
$$\leq c_2 + c_3 z_p^2 \mathrm{Lip}_\sigma^2 \left(\frac{p}{p-1}\right)^2 \exp(\lambda_0^2 t(p-1))$$
$$\times \int_0^t \int_{\mathbf{R}} |p(t-s,x,y)|^2 [\mathrm{E}|u(y,s)|^p]^{2/p} \, dyds.$$

Let $F_p(t) := \sup_{x \in \mathbf{R}}[\mathrm{E}|u(t,x)|^p]^{2/p}$, then it follows that

$$F_p(t) \leq c_2 + c_3 z_p^2 \mathrm{Lip}_\sigma^2 \left(\frac{p}{p-1}\right)^2 \exp(\lambda_0^2 t_0(p-1))$$
$$\times \int_0^t ds F_p(s) \int_{\mathbf{R}} p^2(t-s,x,y) dy$$
$$= c_2 + c_3 z_p^2 \mathrm{Lip}_\sigma^2 \left(\frac{p}{p-1}\right)^2 \exp(\lambda_0^2 t_0(p-1))$$
$$\times \int_0^t ds F_p(s) p(2(t-s),0)$$
$$\leq c_2 + c_4 z_p^2 \mathrm{Lip}_\sigma^2 \left(\frac{p}{p-1}\right)^2 \exp(\lambda_0^2 t_0(p-1))$$
$$\times \int_0^t (t-s)^{\frac{\alpha-1}{\alpha}-1} F_p(s) ds$$

and the result follows by proposition 3.4.

# 6 Conclusion

For the class of equation considered in [1,17], it is known that as time $t$ goes to infinity, the second moment of the solution $E|u(x,t)|^2$ grows like $\exp(constant \times t)$ whenever the initial condition $u_0(x)$ is bounded below; in the same light, we were able to prove the energy growth (second moment growth) of the solution to our class of equation for large time $t$ only in law (in distribution) for the case of a general $L^2$ generator of a Lévy process and say that our solution has a weak (exponential) energy growth for large $t$. For the case of $-(-\Delta)^{\alpha/2}$, the generator of an isotropic stable process, not only did we estimate the $p$ th energy growth bound for all $p > 1$, but were able to give precise exponent of the exponential growth bound in terms of $\alpha$. We also established its existence and uniqueness for all $p > 1$. A further work is to investigate the moment growth property for the gamma-type space-time harmonic function which is also a martingale.

## Acknowledgement

Our appreciation goes to the three anonymous reviewers whose comments helped to improve the content and readability of this article.

## Competing Interests

Author has declared that no competing interests exist.